\documentclass[11pt,A4]{article}

\usepackage{graphicx}
\usepackage{amssymb}
\usepackage{latexsym}

\newtheorem{teo}{Theorem}[section]
\newtheorem{lem}{Lemma}[section]
\newtheorem{defi}{Definition}[section]

\newtheorem{claim}{Claim}[section]

\newenvironment{dem}
{\noindent{\bf Proof:}}
{\hfill$\Box$\\
\par}

\begin{document}

\title{Recognizing [h,2,1] graphs}

\author{Liliana Alc\'on \,\,\, Marisa Gutierrez \thanks{CONICET} \,\,\,
Mar\'\i a P\'\i a Mazzoleni \thanks{Supported by CONICET} \\ \mbox{Departamento de Matem\'atica}\\
Universidad
Nacional de La Plata\\ C. C. 172, (1900) La Plata, Argentina\\
{\tt liliana,marisa,pia@mate.unlp.edu.ar}}


\date{}
\maketitle

\begin{abstract}
An $(h,s,t)$-representation of a graph $G$ consists of a
collection of subtrees of a tree $T$, where each subtree
corresponds to a vertex of $G$ such that (i) the maximum degree of
$T$ is at most $h$, (ii) every subtree has maximum degree at mots
$s$, (iii) there is an edge between two vertices in the graph $G$
if and only if the corresponding subtrees have at least $t$
vertices in common in $T$. The class of graphs that have an
$(h,s,t)$-representation is denoted $[h,s,t]$.\

An undirected graph $G$ is called a $VPT$ graph if it is the
vertex intersection graph of a family of paths in a tree. In this
paper we characterize $[h,2,1]$ graphs using chromatic number. We
show that the problem of deciding whether a given $VPT$ graph
belongs to $[h,2,1]$ is NP-complete, while the problem of deciding
whether the graph belongs to $[h,2,1]-[h-1,2,1]$ is NP-hard. Both
problems remain hard even when restricted to $Split \cap VPT$.
Additionally, we present a non-trivial subclass of $Split \cap
VPT$ in which these problems are polynomial time solvable.
\end{abstract}
\begin{flushleft}
\end{flushleft}

\section{Introduction}
The intersection graph of a set family is a graph whose vertices
are the members of the family, and the adjacency between them is
defined by a non-empty intersection of the corresponding sets.
Classical examples are interval graphs and chordal graphs.\

An $\textbf{interval graph}$ is the intersection graph of a family
of closed intervals on the real line, or equivalently the
intersection graph of a family of subpaths of a path. A
$\textbf{chordal graph}$ is a graph without induced cycles of
length at least four. Gravril \cite{B1} proved that a graph is
chordal if and only if it is the intersection graph of a family of
subtrees of a tree, considering vertex intersection. Both classes
has been widely studied \cite{classes}.

In order to allow larger families of graphs to be represented by
subtrees, several graph classes are defined imposing conditions on
trees, subtrees and intersection sizes  \cite{Jami,B6}. An
$\textbf{(h,s,t)-representation}$ of a graph $G$ consists of a
collection of subtrees of a tree $T$, each subtree corresponding
to a vertex of $G$, such that (i) the maximum degree of $T$ is at
most $h$, (ii) every subtree has maximum degree at mots $s$, (iii)
there is an edge between two vertices in the graph $G$ if and only
if the corresponding subtrees have at least $t$ vertices in common
in $T$. The class of graphs that have an $(h,s,t)$-representation
is denoted $\textbf{[h,s,t]}$. When there is no restriction on the
degree of $T$ or on the degree of the subtrees, we use $h=\infty$
and $s=\infty$ respectively. Notice that $[\infty,\infty,1]$ is
the class of chordal graphs; $[2,2,1]$ is the class of  interval
graphs; $[\infty,2,1]$ and $[\infty,2,2]$ are the well known $VPT$
and $EPT$ graphs \cite{B5}.\

In \cite{Eaton}, the minimum $t$ such that a given graph belongs
to $[3,3,t]$ is studied. In \cite{weakly}, $[4,4,2]$ graphs are
characterized and a polynomial time algorithm for their
recognition is given. In \cite{B7}, the class $[4,2,2]$ is
studied. In \cite{k-edge}, different aspects of $[\infty, 2, t]$
graphs are considered. The relation between the different classes
is analyzed in \cite{B10}. In \cite{B4}, it is shown that the
problem of recognizing $VPT$ graphs is polynomial times solvable,
but the recognition of $EPT$ graphs is an NP-complete problem.\

In this work we focuses in the classes $[h,2,1]$, all them are
subclasses of $VPT$. The problem is deciding  whether a given
$VPT$ graph can be represented as intersection of paths in a tree
with maximum degree $h$. Since $[2,2,1]=Interval$ and
$[3,2,1]=[3,2,2]=EPT \cap chordal$ \cite{B4}, we consider $h\geq
4$. We characterize $[h,2,1]$ graphs using chromatic number. We
show that the problem of deciding whether a given $VPT$ graph
belongs to $[h,2,1]$ is NP-complete, while the problem of deciding
whether the graph belongs to $[h,2,1]-[h-1,2,1]$ is NP-hard. Both
problems remain hard even when restricted to $Split \cap VPT$.
Additionally, we present a non-trivial subclass of $Split \cap
VPT$ in which these problems are polynomial time solvable. In
Section \ref{s:preliminares}, we provide basics definitions and
known results. In Section \ref{s:caract}, we characterize
$[h,2,1]$ graphs for $h\geq 3$. In Section \ref{s:complejidad}, we
present the results about time complexity. Finally, in Section
\ref{s:conclusion} we present some open questions.

\section{Preliminaries}
\label{s:preliminares}
In this paper, all graphs are connected,
finite, simple and loopless. Let $G$ be a graph  with vertex set
$V(G)$ and edge set $E(G)$, the \textbf{open neighborhood}
$\mathbf{N_G(v)}$ of a vertex $v$ is the set of all  vertices
adjacent to $v$. The \textbf{closed neighborhood}
$\mathbf{N_G[v]}$ is $N_G(v)\cup \{v\}$. The \textbf{degree} of
$v$, denoted by $\mathbf{d_G(v)}$, is the cardinality of $N_G(v)$.
For simplicity, when no confusion arise, we omit the subindex $G$
and simply write $N(v)$, $N[v]$ or $d(v)$.

For $S\subseteq V(G)$, $\mathbf{G[S]}$ is the subgraph of $G$
induced by $S$;
 $\mathbf{G-S}$ is a shorthand for $G[V(G)-S]$; and $\mathbf{G-v}$ is used for
 $G-\{v\}$ .

A \textbf{complete set} is a subset of vertices inducing a
complete subgraph. A \textbf{clique} is a maximal complete set.\
The set of cliques of $G$ is denoted by $\mathbf{\mathcal{C}(G)}$.
A $\textbf{stable set}$ is a subset of vertices pairwise
non-adjacent.

The graph $G$ is \textbf{split} if $V(G)$ can be partitioned into
a stable set $S$ and a clique $K$ \cite{classes}. The pair
$\mathbf{(S,K)}$ is  the \textbf{split partition} of $G$. The
vertices in $S$ are called \textbf{stable vertices}, and $K$ is
called the \textbf{central clique} of $G$.  We say that a stable
vertex $s\in S$ is $\textbf{dominated}$ if there exists $s'\in S$
such that $N(s)\subseteq N(s')$. Notice that if $G$ is split then
$\mathcal{C}(G)=\{K, N[s] \mbox{ for } s\in S\}$.

A \textbf{VPT-representation} of $G$, denoted by $\langle
\mathcal{P}, T\rangle$,  is an $(\infty,2,1)$- representation.
This means that $\mathcal{P}$ is a family $(P_v)_{v \in V(G)}$ of
subpaths of a host tree $T$ satisfying that two vertices $v$ and
$v'$ of $G$ are adjacent if and only if $P_v$ and $P_{v'}$ have at
least one vertex in common. If $q$ is a vertex of the host tree
$T$, then $\mathbf{P[q]}$ denote the set $\{P\in \mathcal{P} /
q\in V(P)\}$ and $\mathbf{C_q}$ denote the complete set $\{v\in
V(G) / q\in V(P_v) \}$. Notice that for every clique  $C$ of $G$,
there exists $q \in V(T)$ such that $C=C_q$.

\begin{defi}\cite{B4} Let $C\in \mathcal{C}( G)$. The $\textbf{branch
graph}$  of $G$ for the clique $C$ denoted by $\mathbf{B(G/C)}$ is
defined as follows: the vertex set $V(B(G/C))$ contains the
vertices of $V(G)\setminus C$ adjacent to some vertex of $C$. Two
vertices $v$ and $w$  are adjacent in   $B(G/C)$ if and only if
\begin{enumerate}
\item $vw\notin E(G)$;

\item there exists a vertex of $C$ adjacent to both $v$ and $w$;
and

\item there exist vertices $v'$ and $w'$ of $C$ such that $v'$ is
adjacent to $v$ and non-adjacent to $w$, and $w'$ is adjacent to
$w$ and non-adjacent to $v$.
\end{enumerate}
\end{defi}

Let $q\in V(T)$, with $N_{T}(q)=\{y_{1},y_{2},..,y_{h}\}$. We call
$\textbf{branches of T at q}$ to the connected components of
$V(T)$ - $\{q\}$. Observe that each $y_{i}$ is contained in a
different branch which will be called $T_{i}$.

The graph $G$ is $\textbf{k-colorable}$  if its vertices can be
colored with at most $k$ colors in such a way that  no two
adjacent vertices share the same color. The $\textbf{chromatic
number}$ of $G$, denoted by $\chi(G)$, is the smallest number of
colors needed to coloring $G$.

\begin{teo}\cite{karp}
Let $G$ be a graph and $k\geq 3$. Deciding whether $G$ is
$k$-colorable is an NP-complete problem.
\end{teo}

A graph $G$ is \textbf{perfect} if and only if $G$ is
$\{$$C_{2n+1}$, $\bar{C}_{2n+1}$, with $n\geq 2$ $\}$- free
\cite{B0}.

\begin{teo}\cite{gro}\label{t:per} Let $G$ be a perfect
graph and $k\geq 3$. Deciding whether $G$ is $k$-colorable is a
polynomial time solvable problem.
\end{teo}

\section{Characterization of [h,2,1], for $h\geq 3$}
\label{s:caract}

In this section we present a characterization of the $VPT$ graphs
that can be represented in a tree with maximum degree at most $h$.
The characterization is given  in terms of the chromatic number of
the branch graphs. The following three lemmas are  fundamental
tools in the proof of the main Theorems \ref{t:coloreable} and
\ref{t:cromatic}.

\begin{lem}
\label{l:branch} Let $\langle \mathcal{P}, T \rangle$ be a $VPT$
representation of $G$. Let $C\in \mathcal{C}(G)$ and $q\in V(T)$
such that $C=C_{q}$.  If $v\in V(B(G/C))$ then $P_v$ is contained
in some branch of  $T$ at $q$.  If $v$ is adjacent to $w$ in
$B(G/C)$ then $P_v$ and $P_{w}$ are not contained in a same branch
of $T$ at $q$.
\end{lem}

\begin{dem}
If $v\in V(B(G/C))$ then $v\notin C$. It follows that $q\notin
V(P_{v})$, thus $P_v$ is contained in some branch $T_i$ of $T$ at
$q$. Let $w\in N_{B(G/C)}(v)$ and assume for a contradiction that
$P_{v}$ and $P_w$ are contained in the same branch $T_{i}$. Let
$x$ and $y$ be the vertices of $P_v$ and $P_w$ respectively  at
minimum distance from $q$. Since there exists a vertex of $C$
adjacent to $v$ and $w$, there exists  a path in $T$ containing
$q$, $x$ and $y$. We can assume, without loss of generality, that
$x$ is between $q$ and $y$ or that $x=y$. In both cases, $N(w)\cap
C \subseteq N(v)\cap C$. This contradicts the fact that $v$ and
$w$ are adjacent in the branch graph.
\end{dem}

\begin{lem}
\label{l:grado} Let $\langle \mathcal{P}, T \rangle$ be a $VPT$
representation of $G$. Let $C\in \mathcal{C}(G)$ and $q\in V(T)$
such that $C=C_{q}$. If $d_{T}(q)=h$, then $B(G/C)$ is
$h$-colorable.
\end{lem}

\begin{dem}
 Let $T_{1}, T_{2},..,T_{h}$ be the branches
of T at q. By Lemma \ref{l:branch}, if we color each vertex $v$ of
$B(G/C)$ with the index $i$ of the branch $T_i$ containing $P_v$,
then we obtain a proper coloring of $B(G/C)$. Since there are $h$
branches, $B(G/C)$ is $h$-colorable.
\end{dem}

\begin{lem}
\label{l:repre} Let $\langle\mathcal{P},T\rangle$ be a $VPT$
representation of $G$. Consider $q\in V(T)$ with $d_{T}(q)=h\geq
4$. Assume  there exist $y_1, y_2\in N_T(q)$ such that for all
$v\in V(G)$, $\{y_{1}, y_{2}\} \nsubseteq V(P_{v})$. Then there
exists a $VPT$ representation $\langle \mathcal{P}',T'\rangle$ of
$G$ with $V(T')=V(T)\cup \{a_{q}\}$, $a_q \notin V(T)$, and

$$d_{T'}(x) =\left\{
\begin{array}{ll}
3,  \ \ \ &\mbox{if\ }x = a_{q}\\
 \
h-1, \ \ \ &\mbox{if\ }x = q\\
 \
d_{T}(x), \ \ \ &\mbox{if\ }x \in V(T')\setminus \{q,a_{q}\}.
\end{array}
\right.$$
\end{lem}

\begin{figure}[h]
\centering{
\includegraphics[height=1.6in,width=3in]{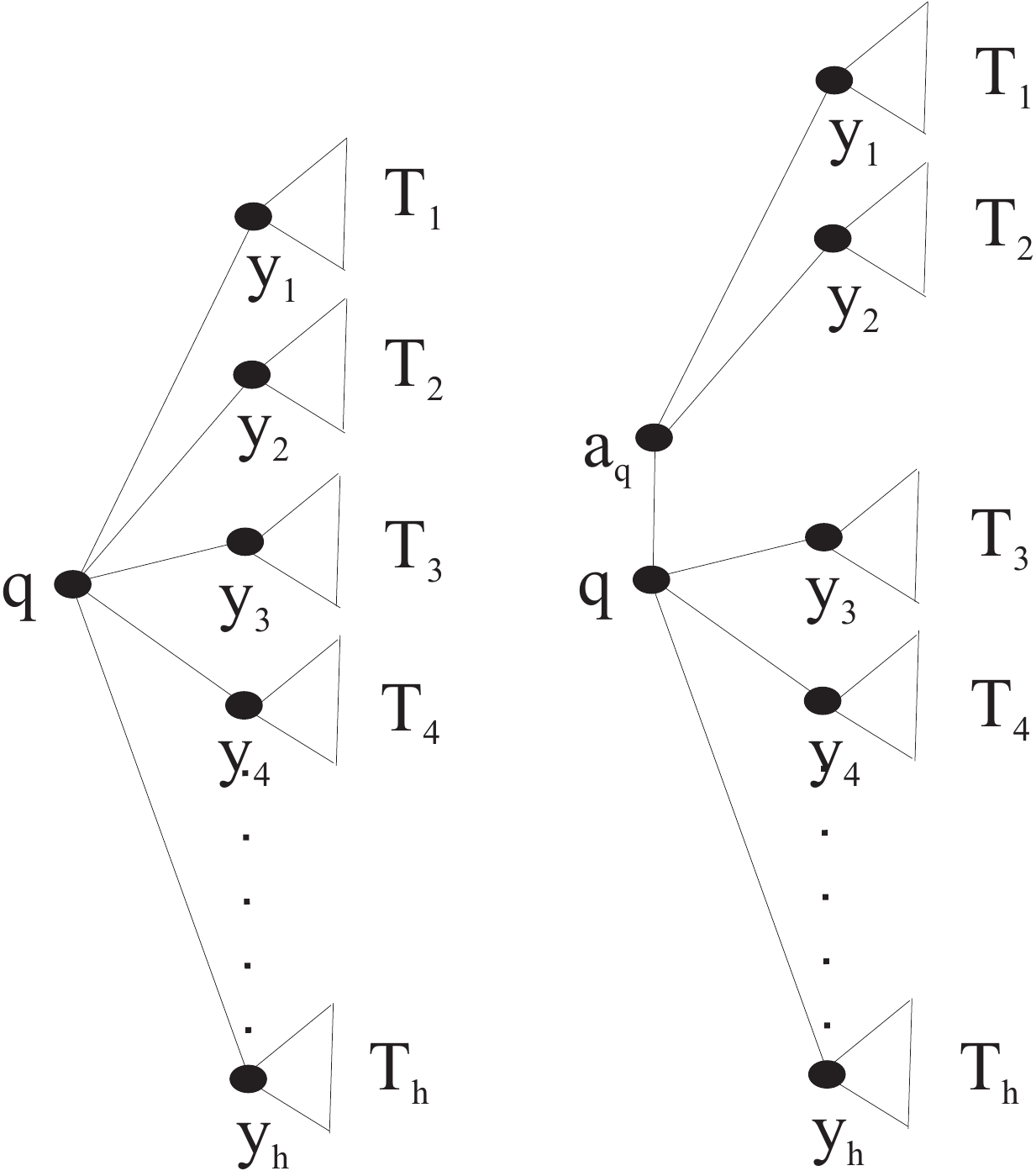}}
\caption{The  degree of $q$ in the tree $T$ on the left is $h$.
The degree of $q$ in the tree $T'$ on the right is $h-1$.}
\label{f:TyT'}
\end{figure}

\begin{dem}
We obtain the  $\langle \mathcal{P}',T'\rangle$ representation of
$G$ as follows (Please refer to Figure \ref{f:TyT'}): the set of
vertices of
 $T'$ is $V(T) \cup
\{a_{q}\}$, where $a_{q}$ is a new vertex not in $V(T)$. The set
of edges is $(E(T) \setminus \{ y_{1}q, y_{2}q \}) \cup \{
y_{1}a_{q}, y_{2}a_{q}, qa_{q} \}$. Observe that the degree of
each vertex  $x\in V(T')$ is the required in the statement of the
present lemma.

Now we define the  paths $P'_{v}$ for $v\in V(G)$: if $y_{1}$ and
$q$ or $y_{2}$ and $q$ belong to $V(P_{v})$ then
$V(P'_{v})=V(P_{v})\cup \{a_{q}\}$. In any other case,
$V(P'_{v})=V(P_{v})$. Since $\{y_{1},q,y_{2}\} \subsetneqq
V(P_{v})$, we have that each $V(P'_{v})$ induces a path in $T'$.
Moreover, since all the paths where vertex $a_{q}$ was added had
vertex $q$ in common, it is clear that, for any pair of vertices
$v, w\in V(G)$, $V(P_{v})\cap V(P_{w})\neq \emptyset$ if and only
if $V(P'_{v})\cap V(P'_{w})\neq \emptyset$. It  follows that
$\langle \mathcal{P}',T'\rangle$ is a $VPT$ representation of $G$
and the implication is proven.
\end{dem}

\begin{teo}
\label{t:coloreable} Let $G\in VPT$ and  $h\geq 3$. The graph $G$
belongs to  $[h,2,1]$ if and only if  $B(G/C)$ is $h$-colorable
for every $C\in \mathcal{C}(G)$. The direct implication is true
also for $h=2$.
\end{teo}

\begin{figure}[h]
 \centering{
  \includegraphics[height=2.2in,width=3in]{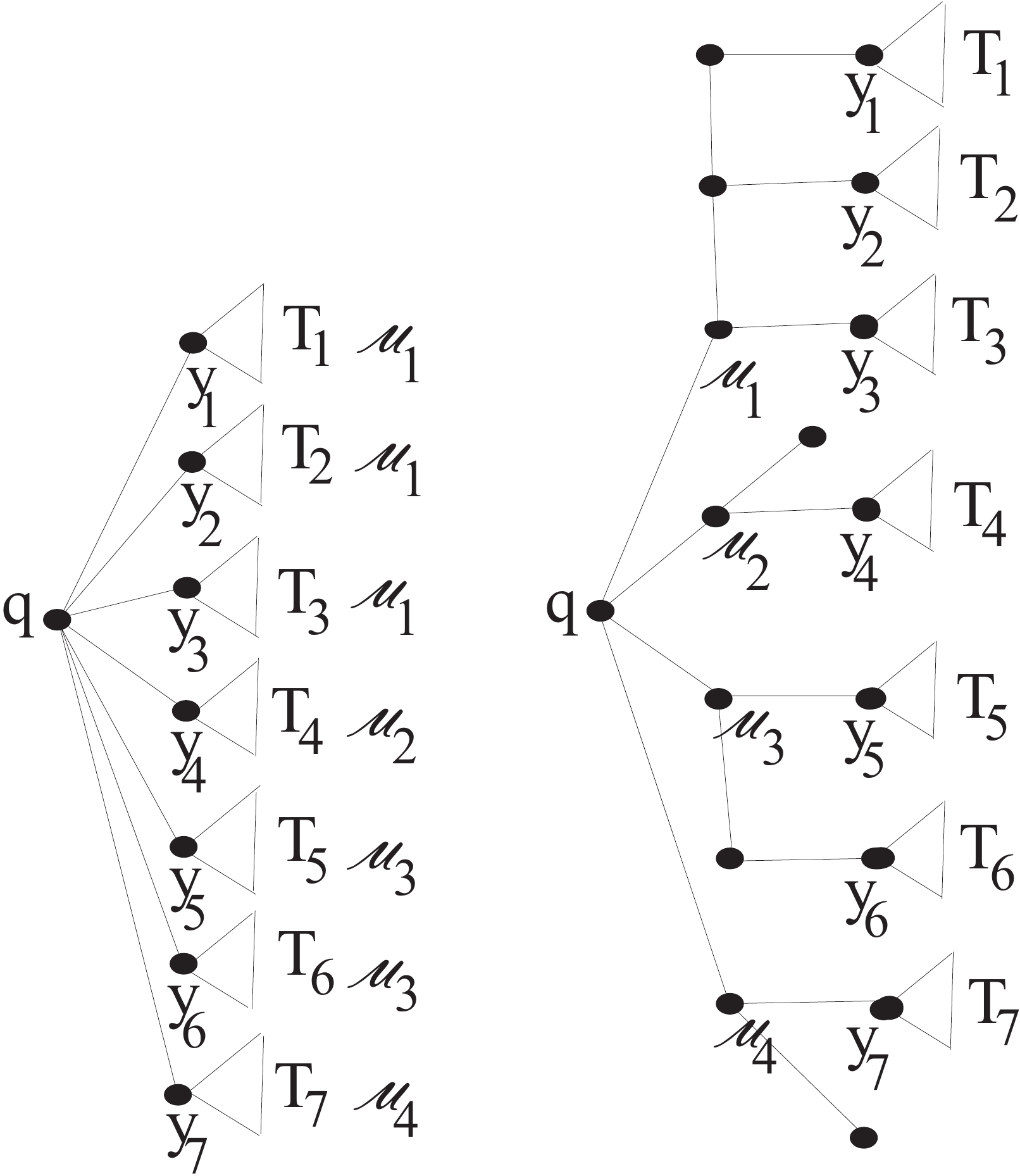}}
  \caption{$d_{T}(q)=7$ and $B(G/C_{q})$ is $4$-colorable.}
  \label{f:mus}
  \end{figure}

\begin{dem}
Let $\langle \mathcal{P},T \rangle$ be an $(h,2,1)$-representation
of $G$ with $h\geq 2$. Assume $C\in \mathcal{C}(G)$, then there
exists $q\in V(T)$ such that $C=C_{q}$. Since $d_{T}(q)\leq h$, by
Lemma \ref{l:grado}, $B(G/C)$ is $h$-colorable.

The reciprocal  implication for $h=3$ was proven  by Golumbic and
Jamison in \cite{B4}; then we  assume $h\geq 4$.

 Let $\langle \mathcal{P},
T\rangle$ be a $VPT$ representation of $G$. We will prove that $G$
admits an $(h,2,1)$-representation.

 We proceed by induction on the number $k$ of vertices
of $T$ whose degree exceeds $h$. If $k=0$ we are done.

If $k> 0$,  there exists a vertex $q$ of $T$ with degree $d>h$.
Say $N_{T}(q)=\{y_{1},y_{2},...,y_{d}\}$ and for every $i$, $1\leq
i\leq d$, let $T_{i}$ be the branch of $T$ at $q$ containing the
vertex $y_{i}$.

If by repeatedly applying the Lemma \ref{l:repre} we can obtain a
$VPT$ representation $\langle P',T'\rangle$ of $G$ with
$d_{T'}(q)<h$, then the implication is proven by induction  since
no vertex of $T'$ increases its degree.

In other case, we can assume that for any pair of vertices
$y_{i}$, $y_{i'}$ belonging to $N_{T}(q)$, there exists at least
one $v\in V(G)$ such that $\{y_{i}, y_{i'}\}\subseteq V(P_{v})$.

Notice that this implies that  $C_{q}=$$\{v\in V(G) / q\in
V(P_v)\}$ is a clique of $G$.

We will consider two cases.

Case 1: for every $i$, $1 \leq i\leq d$, there exists $v_i \in
V(G)$ such that $P_{v_{i}}$ is totally contained in the branch
$T_i$ and $y_i \in V(P_{v_{i}})$. Observe that each $v_{i}$ must
be a vertex of $B(G/C_{q})$. Since $B(G/C_{q})$ is $h$-colorable,
we can partitioned the set $\{y_{1},y_{2},...,y_{d}\}$ in $h$
subsets $D_j$, $1\leq j \leq h$, each one containing the vertices
$y_i$ for which the associated vertex $v_i$ has color $j$.

We  obtain a new $VPT$ representation $\langle P',T'\rangle$ of
$G$ as follows. The tree $T'$ is obtained from $T$ by means of the
following procedure (in Figure \ref{f:mus} we offer an example):
1) remove the edges $qy_{i}$, $1 \leq i\leq d$; 2) add $h$ new
vertices $\mu_{j}$, $1\leq j \leq h$; 3) add the edges $q\mu_{j}$,
$1\leq j \leq h$; and finally, to connect the vertices $\mu_{j}$
with the vertices $y_i$, 4) add for every $j$, $1 \leq j \leq h$,
a binary tree rooted at the vertex $\mu_{j}$ and  with the
vertices of $D_j$ as leaves. The rest of the tree T remains
unchanged.

The only paths of $\mathcal{P}$ which are modified  to obtain the
paths of $\mathcal{P'}$  are those $Q\in P[q]$. If $Q$ has $q$ as
an endpoint, then we obtain $Q'$ by replacing in $Q$ the edge
$qy_{i}$ by the unique subpath of $T'$ linking $q$ and $y_{i}$. If
$Q$ has $q$ as an internal vertex, then there exist $i$ and $i'$
such that $Q$ contains the edges $qy_{i}$ and $qy_{i'}$. Notice
that the existence of $Q$  implies that  $v_i$ and $v_{i'}$ are
adjacent in $B(G/C_{q})$; thus they have different colors, say $j$
and $j'$. Therefore, we obtain $Q'$  by replacing in $Q$ the edges
$qy_{i}$  and  $qy_{i'}$ by the only  subpath of $T'$ linking
$y_{i}$, $q$ and  $y_{i'}$.

It is easy to see that this construction leaves the intersection
graph of paths unchanged while reducing the number of tree
vertices of degree greater than h. So, by induction, the
implication is proven.\

Case 2: there exists $i$, $1 \leq i\leq d$,  such that every path
$P \in \mathcal{P}$ containing $y_i$ is not contained in the
branch $T_i$. Thus, every path $P \in \mathcal{P}$ containing
$y_i$  contains also $q$. Therefore, we can contract the edge
$qy_{i}$ to obtain a new $VPT$ representation of $G$ and repeat
this as many times as needed to get a representation which is in
Case 1. Notice that in this procedure some vertices of $T$
disappear, and that the degree of $q$ may  increase, but the
number of vertices  whose degree exceeds $k$ does not grow, thus
the proof follows by induction as in the previous case.
\end{dem}

Observe that the reciprocal implication of Theorem
\ref{t:coloreable} is false for $h=2$; consider,  by instance, the
graph $T^{3}_{2}$ which consists of one central vertex and 3 edge
disjoint paths of 2 edges each intersecting only on the central
vertex. It is easy to see that $T^{3}_{2}\in VPT$ and $B(G/C)$ is
$2$-colorable  for all $C\in \mathcal{C}(G)$, but $T^{3}_{2}\notin
[2,2,1]$.

\begin{teo}
\label{t:cromatic}
 Let $G\in VPT$ and  $h\geq 4$. The graph $G$ belongs to
 $[h,2,1]-[h-1,2,1]$ if and only if $Max_{C\in
\mathcal{C}(G)}(\chi(B(G/C)))=h$. The reciprocal implication is
also true for $h=3$.
\end{teo}

\begin{dem}
By Theorem \ref{t:coloreable}, $G\in [h,2,1]$ if and only if
$Max_{C\in \mathcal{C}(G)}(\chi(B(G/C)))\leq h$. On the other
hand, by the same Theorem \ref{t:coloreable},  $G \notin
[h-1,2,1]$ if and only if $Max_{C\in
\mathcal{C}(G)}(\chi(B(G/C)))> h-1$. Therefore, the proof follows.
\end{dem}

\section{Complexity}
\label{s:complejidad} In this Section we  prove that the problem
of deciding whether a given graph belongs to $[h,2,1]$  for $h
\geq 3$ is NP-complete. We also show that recognizing
$[h,2,1]-[h-1,2,1]$ for $h\geq 4$ is NP-hard.  Our results prove
that both problems remain  difficult even when restricted to the
class $VPT \cap Split$ without dominated stable vertices.

First we state the following fundamental property of $VPT \cap
Split$ graphs which is used in the proof of Theorems \ref{t:hard}
and \ref{t:complete}.

\begin{lem}
\label{l:1-col}Let $s$ be a stable vertex of a $VPT \cap Split$ graph
$G$. The branch graph $B(G/N[s])$ is $1$-colorable.
\end{lem}

\begin{dem}
Let $\langle P,T \rangle$ be a $VPT$ representation of $G$ such
that $P_s$ is a one vertex path in a leaf $y$ of $T$, in other
words $V(P_s)=\{y\}$ where $y$ is a leaf of $T$. Thus $N[s]$ is
the clique $C_y$. Since $d_T(y)=1$, by Lemma \ref{l:grado},
$B(G/N[s])$ is $1$-colorable.
\end{dem}

For the NP-completeness proof, we use a reduction from the
chromatic number problem \cite{karp}.

Given a graph $G$ we construct in  polynomial time
 a graph $\widehat{G}\in VPT\cap Split$ without dominated stable vertices,
 in such a way that $\chi(G)=h$
 if and only if $\widehat{G} \in [h,2,1]-[h-1,2,1]$.

Let $V(G)=\{v_{1},v_{2},...,v_{n}\}$,  we define the graph
$\widehat{G}$ by means of its $VPT$ representation $\langle
\mathcal{P},T \rangle$ as follows: the tree $T$ is a star with
central vertex $q$ and  leaves $y_i$ for $1\leq i \leq n$.\

The path family $\mathcal{P}$ contains: a one vertex path $P_{i}$
with $V(P_{i})= \{y_{i}\}$, for each $1\leq i\leq n$; a path
$P_{ij}$ with $V(P_{ij})=\{y_{i}, q, y_{j}\}$, for each $1\leq
i<j\leq n$ such that $v_{i}v_{j}\in E(G)$; a path $P_{iq}$ with
$V(P_{iq})=\{q, y_{i}\}$, for each $1\leq i\leq n$ such that
$d_{G}(v_{i})=1$.\

We call each vertex of $\widehat{G}$ as the corresponding path of
$\mathcal{P}$.\

\begin{figure}[h]
\centering{
\includegraphics[height=2in,width=3in]{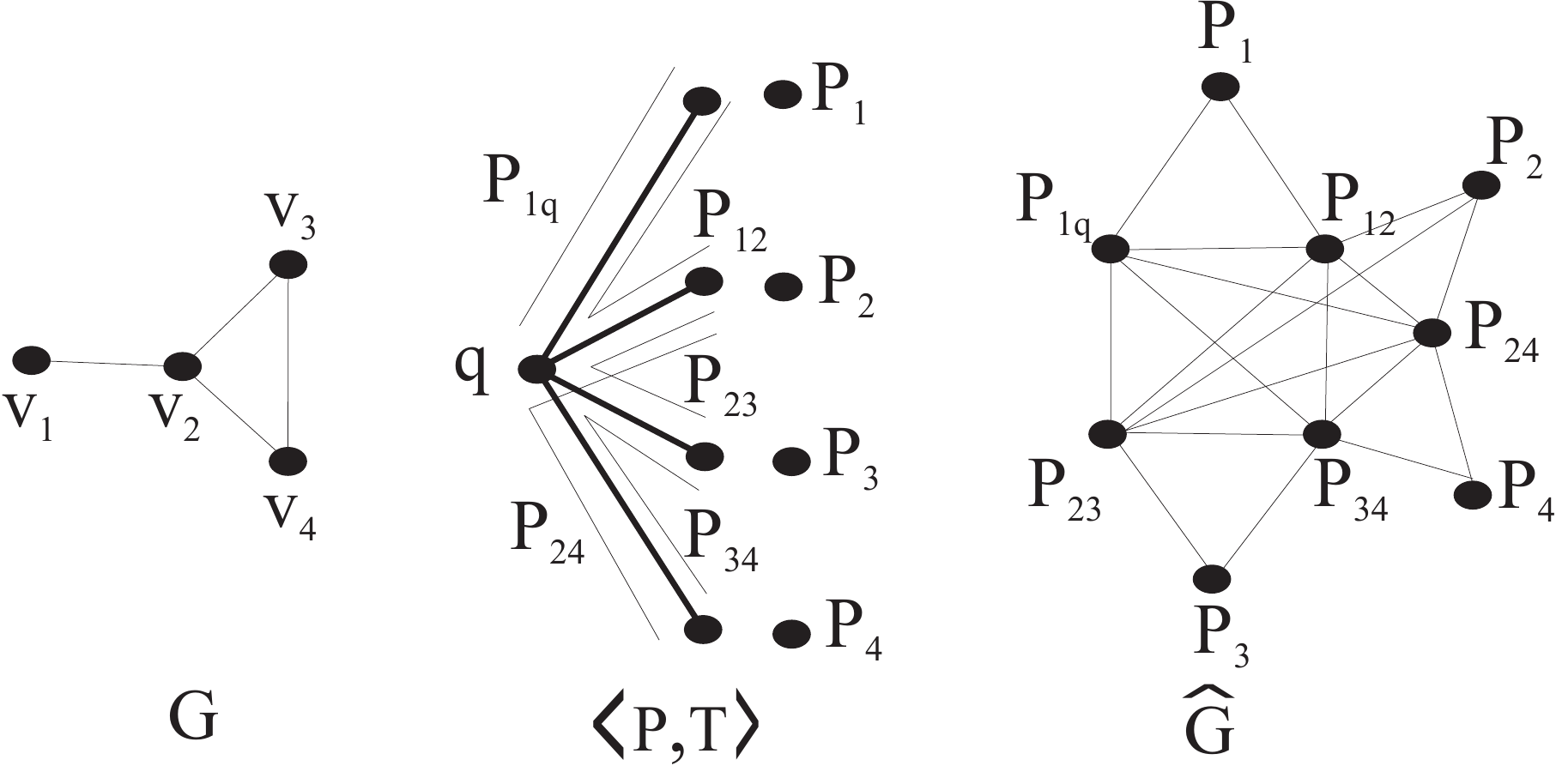}}
\caption{A graph $G$, the $VPT$ representation of $\widehat{G}$
and the graph $\widehat{G}$.}
\label{f:ejemG}
\end{figure}

In Figure \ref{f:ejemG} we offer an example of a graph $G$, the
$VPT$ representation of $\widehat{G}$ and the graph $\widehat{G}$
obtained.

\

Notice that $\widehat{G}$ is a split graph with the vertex set
partitioned in a stable set of size $n=\mid V(G)\mid$
corresponding to the one vertex paths $P_i$; and a  central clique
of size $\mid E(G)\mid + \mid \{v \in V(G) / d_G(v)=1\}\mid$
corresponding to the remaining paths, all of which contain  the
vertex $q$ of $T$, thus this clique is $C_q$. The other cliques of
$\widehat{G}$ are the cliques $C_{y_i}$ for $1 \leq i \leq n$ each
one corresponding to the paths containing the vertex $y_i$ of $T$
respectively. The graph $\widehat{G}$ has no more cliques. In
addition, every stable vertex $P_i$ of $\widehat{G}$ is
non-dominated.

Observe that the branch graphs $B(\widehat{G}/C_{y_i})$ are
described in Lemma \ref{l:1-col}, the following claim does for
$B(\widehat{G}/C_{q})$.

\begin{claim}
\label{c:1} If $\widehat{G}$ is the graph obtained from $G$ as
above, then $B(\widehat{G}/C_{q})=G$.
\end{claim}

\begin{dem}
Notice that $B(\widehat{G}/C_{q})$ has exactly $n$ vertices: $P_i$
for $1\leq i \leq n$.

We will see that $P_i$ and $P_j$ are adjacent in
$B(\widehat{G}/C_{q})$ if and only if $v_i$ and $v_j$ are adjacent
in $G$. If  $P_iP_{j}\in E(B(\widehat{G}/C_{q}))$ then  there
exists a vertex of $C_{q}$ adjacent to both $P_{i}$ and $P_{j}$.
Then, there exists a path $P_{ij}\in \mathcal{P}$, thus
$v_{i}v_{j}\in E(G)$. Reciprocally, let $v_{i}v_{j}\in E(G)$.
Notice that $P_i$ and $P_j$ are non-adjacent in $\widehat{G}$; and
$P_{ij}$ is a vertex of $C_q$ adjacent to $P_i$ and to $P_j$ in
$\widehat{G}$. Let us see that there exists a vertex of $C_q$
adjacent to $P_i$ and non-adjacent to $P_j$. Indeed, if
$d_G(v_i)=1$ then  the wanted vertex of $C_q$ is $P_{iq}$. If
$d_G(v_i)>1$ then $v_i$ must have a neighbor $v_l$ with $l\neq j$,
thus  the wanted  vertex of $C_q$ is $P_{il}$. In an analogous
way, there exists a vertex of $C_q$ adjacent to $P_j$ and
non-adjacent to $P_i$. We have proved that $P_i$ and $P_j$ are
adjacent in $B(\widehat{G}/C_{q})$. We conclude that
$B(\widehat{G}/C_{q})=G$.
\end{dem}

The reduction from chromatic number is complete using the next
claim.

\begin{claim}
\label{c:2}
 Let $\widehat{G}$ be the graph obtained from $G$ as
above and  $h \geq 4$. The graph   $\widehat{G}$ belongs to
$[h,2,1]-[h-1,2,1]$ if and only if
 $\chi(G)=h$.
 \end{claim}

\begin{dem}
By  Lemma \ref{l:1-col} and Claim \ref{c:1}, $Max_{C\in
\mathcal{C}(\widehat{G})}\chi(B(\widehat{G}/C))=
\chi(B(\widehat{G}/C_q))=\chi (G)$. Hence, by Theorem
\ref{t:cromatic}, $\widehat{G}$ belongs to $[h,2,1]-[h-1,2,1]$ if
and only if
 $\chi(G)=h$.
\end{dem}

We have proved the following theorem.

\begin{teo}
\label{t:hard} Let  $G\in VPT\cap Split$ without dominated stable
vertices, and $h\geq 4$. Decide whether $G\in [h,2,1]-[h-1,2,1]$
is an NP-hard problem.
\end{teo}

In addition, since an $(h,2,1)$-representation is a polynomial
certificate of belonging to $[h,2,1]$; using Theorem
\ref{t:coloreable} and the construction above, we have proved the
following result.

\begin{teo}
\label{t:complete} Let $G\in VPT\cap Split$  without dominated
stable vertices, and  $h\geq 3$. Decide whether $G\in [h,2,1]$ is
an NP-complete problem.
\end{teo}

We notice that Theorem \ref{t:complete} for $h=3$ has been
previously proved  in \cite{B4}.

\subsection{A polynomial time solvable subclass }

We have proved that deciding whether a given $VPT\cap Split$ graph
without dominated stable vertices admits a representation as
intersection of paths of a tree with maximum degree $h$ is an
NP-complete problem. In what follows we describe a non-trivial
subclass of $VPT \cap Split$ without dominated stable vertices
where the problem is polynomial time solvable.

\begin{figure}[h]
\centering{
\includegraphics[height=1.5in,width=3in]{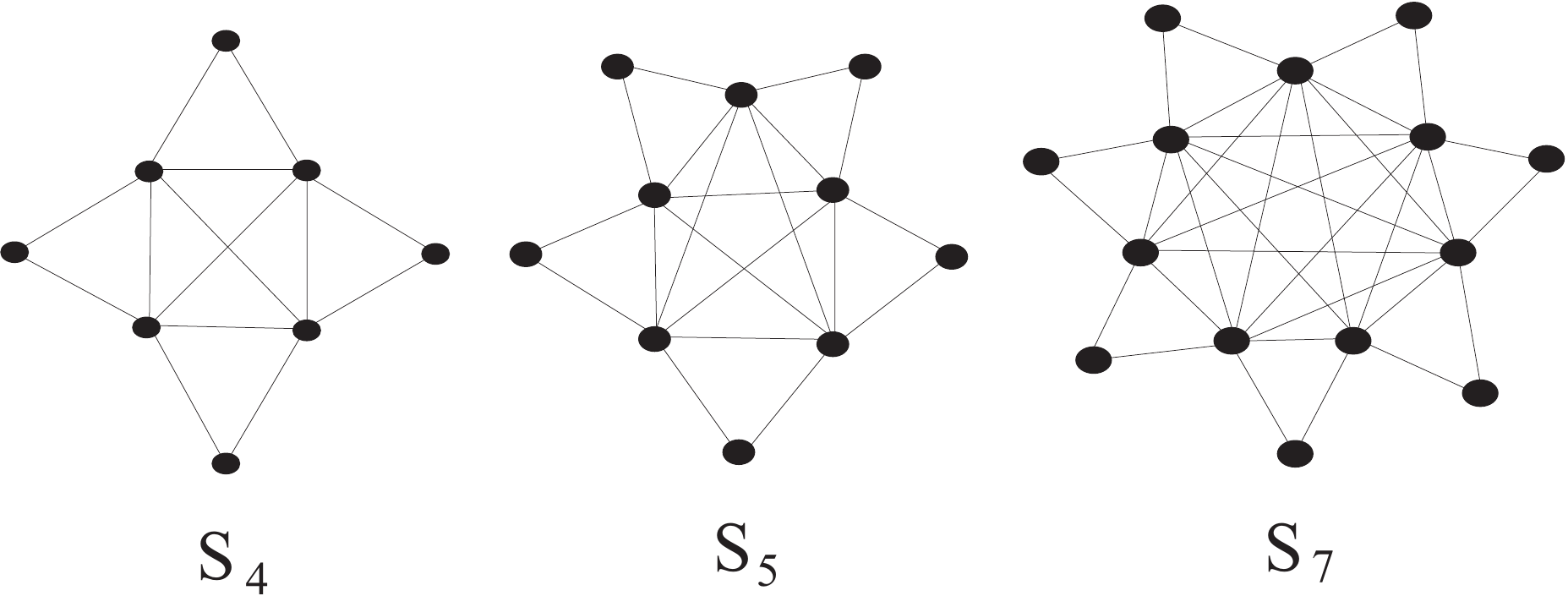}}
\caption{The sun graphs $S_{4}$, $S_{5}$ and $S_{7}$.}
\label{f:estrellas}
\end{figure}

For $n\geq 4$, a $\textbf{n-sun}$, denoted by
\textbf{S}$_\textbf{n}$, is a split graph with stable set
$\{s_{1},s_{2},..,s_{n}\}$, central clique
$\{v_{1},v_{2},..,v_{n}\}$, $N(s_i)=\{v_{i},v_{i+1}\}$ for $1\leq
i\leq n-1$, and $N(s_n)=\{v_n, v_1 \}$. See Figure
\ref{f:estrellas}.\

Let $G$ be a split graph with partition $(S,K)$.
 We say that $G$ belongs to \textbf{$SVS$} (special $VPT$ subclass)
 whenever

 \begin{itemize}
    \item $G \in VPT$,
    \item for all $v\in K$, $|N(v)\cap S|\leq 2$, and
    \item if $S_k$, with $k\in \{4, 2n+1$ for $n\geq 2\}$, is induced in $G$
    then there exists $v\in K$ such that
    $v$ is adjacent to two non-consecutive vertices of the stable set of $S_{k}$.
 \end{itemize}

The class $SVS$ is not trivial, in the sense that it includes
graphs in $[h,2,1]$ for all $h\geq 4$.\

\begin{figure}[h]
\centering{
\includegraphics[height=2.2in,width=3in]{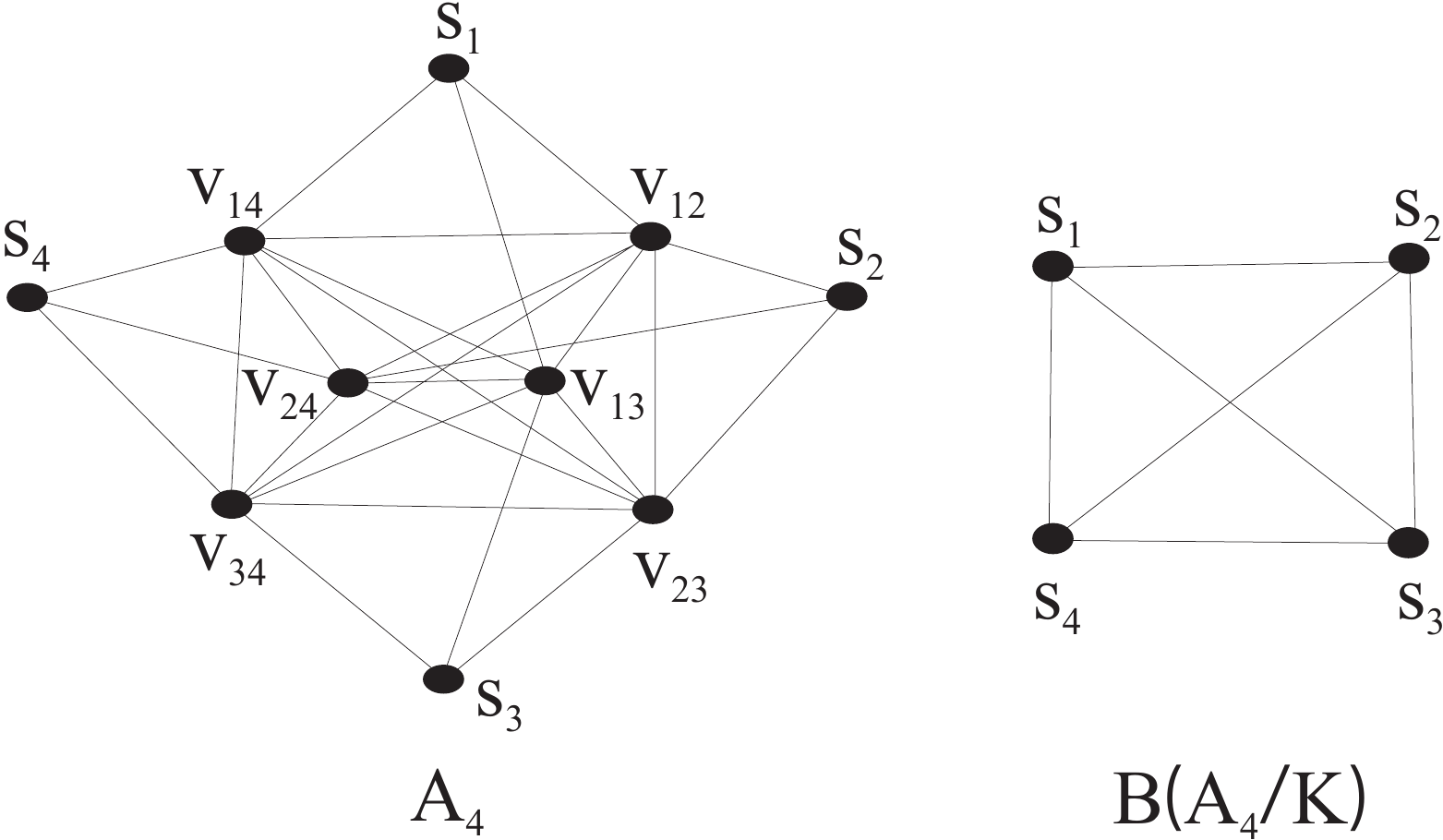}}
\caption{The graph $A_{4}$ belongs to $SVS$ and $A_{4}\in
[4,2,1]-[3,2,1]$.} \label{f:Rn}
\end{figure}

For example, let $n\geq 4$ and let $A_{n}$ (see \cite{B10}) be a
split graph with partition $(S,K)$, where $S=\{s_{1},..,s_{n}\}$,
$K=\{v_{ij}, 1\leq i<j\leq n\}$ and $N(v_{ij})=\{s_{i}, s_{j}\}$,
for all $1\leq i<j\leq n$. It is clear that $A_{n}$ belongs to
$SVS$, and $B(A_{n}/K)=K_{n}$ with $V(K_{n})=\{s_{1},..,s_{n}\}$.
Hence, by Theorem \ref{t:cromatic}, $A_{n}\in [n,2,1]-[n-1,2,1]$.
(As an example see Figure \ref{f:Rn}).

The following two lemmas are used in the proof of the main Theorem
\ref{t:poly} which proves that in the class $SVS$ the graphs
belonging to $[h,2,1]$ can be recognized efficiently.

\begin{lem}
\label{l:S-n}Let $G\in VPT \cap Split$ with partition $(S,K)$ such that
for all $v\in K$, $|N(v)\cap S|\leq 2$, and let $n\geq 4$. If
$B(G/K)$ has an induced $C_n$ then  $G$ has an induced $S_{n}$.
\end{lem}

\begin{dem}
Let $\langle \mathcal{P},T \rangle$ be a $VPT$ representation of
$G$ and $q\in V(T)$ such that $K=C_q$. Let $C_n$ be an induced
cycle of $B(G/K)$ with vertices $s_1$, $s_2$,...,$s_n$. It is
clear that every $s_i \in S$. Since $s_i$ is adjacent to $s_{i+1}$
in  $B(G/K)$,  there exists $v_i \in K$ such that $v_i$ is
adjacent to $s_i$ and to $s_{i+1}$ in $G$. Since, for all $v\in
K$, $|N(v)\cap S|\leq 2$, if $i \neq i'$ then $v_i \neq v_{i'}$,
thus $s_1$, $s_2$,...,$s_n$, $v_1$, $v_2$,...,$v_n$ induce a
$n$-sun in $G$ and the proof is completed.
\end{dem}

\begin{lem}
\label{l:perf} If $G\in SVS$ then every branch graph of $G$ is
perfect.
\end{lem}

\begin{dem}
Let $(S,K)$ be a split partition of $G$. By Lemma \ref{l:1-col},
if $s\in S$ then $B(G/N[s])$ is perfect. Assume for a
contradiction that $B(G/K)$ is not perfect, then $B(G/K)$ has
induced an odd cycle or the complement of an odd cycle. Since the
complement of $C_5$ is $C_5$; and the complement of any odd cycle
of size 7 or more has an induced $C_4$, it follows that $B(G/K)$
has an induced $C_k$, for some $k\in \{4, 2n+1$ for $n\geq 2\}$.
Therefore, by Lemma \ref{l:S-n}, $G$ has an induced $S_k$. Since
$G\in SVS$, there exists $v\in K$ such that $v$ is adjacent to two
non-consecutive vertices $s$ and $s'$ of the stable set of
$S_{k}$. Notice that the existence of $v$ implies that the
vertices $s$ and $s'$ are adjacent in $B(G/K)$. This contradicts
the fact that $C_{k}$ is an induced cycle of $B(G/K)$.\end{dem}

\begin{teo}
\label{t:poly} Let $G\in SVS$ and $h\geq 4$. Decide whether $G$
belongs to $[h,2,1]-[h-1,2,1]$ is polynomial time solvable.
\end{teo}

\begin{dem}
Given $G\in SVS$, in order to determinate if $G\in
[h,2,1]-[h-1,2,1]$, by Theorem \ref{t:coloreable}, it is enough to
calculate the chromatic number of $B(G/K)$, where $K$ is the
central clique of $G$. Notice that the branch graph $B(G/K)$ can
be obtained in polynomial time. On the other hand, by  Lemma
\ref{l:perf}, $B(G/K)$ is perfect. Thus, by Theorem \ref{t:per},
its chromatic number can be calculated in polynomial time.
\end{dem}

\section{Future work}
\label{s:conclusion} In this paper we  give a characterization of
the $[h,2,1]$ graphs, with $h\geq 3$. In addition, we prove that
recognizing this class is NP-complete and show a family, called
$SVS$, in which this problem is polynomial time solvable. We are
working in describing a larger subclass of $VPT$ graphs where this
problem remains polynomial. On the other hand, we are analyzing
the possibility of extending the techniques used in the present
paper to characterize the classes $[h,2,2]$.



\end{document}